\documentclass[12pt]{article}
\usepackage{amsmath}
\usepackage{amssymb}
\usepackage[dvips]{graphicx}
\usepackage{epsfig}
\font\tencmmib=cmmib10 \skewchar\tencmmib '60
\newfam\cmmibfam
\textfont\cmmibfam=\tencmmib

\def\bbox{\quad\hbox{\vrule \vbox{\hrule \vskip2pt \hbox{\hskip2pt
\vbox{\hsize=1pt}\hskip2pt} \vskip2pt\hrule}\vrule}}
\def\lessim{\ \lower4pt\hbox{$
\buildrel{\displaystyle <}\over\sim$}\ }
\def\gessim{\ \lower4pt\hbox{$\buildrel{\displaystyle >}
\over\sim$}\ }

\def\si{\bar{\sigma}}

\def\O{{\cal O}}

\def\la{\langle}
\def\ra{\rangle}

\def\qed{\hfill\break\rightline{$\bbox$}}
\newcommand{\e}{\mathbb{E}}

\newtheorem{lemma}{Lemma}
\newtheorem{theorem}{Theorem}

\makeatletter
\@addtoreset{equation}{section}

\makeatother

\font\tencmmib=cmmib10 \skewchar\tencmmib '60
\newfam\cmmibfam
\textfont\cmmibfam=\tencmmib

\def\bbox{\quad\hbox{\vrule \vbox{\hrule \vskip2pt \hbox{\hskip2pt
\vbox{\hsize=1pt}\hskip2pt} \vskip2pt\hrule}\vrule}}
\def\lessim{\ \lower4pt\hbox{$
\buildrel{\displaystyle <}\over\sim$}\ }
\def\gessim{\ \lower4pt\hbox{$\buildrel{\displaystyle >}
\over\sim$}\ }

\def\go0{\to 0}

\def\la{\langle}

\def\leftitem#1{\item{\hbox to\parindent{\enspace#1\hfill}}}

\def\qed{\hfill\break\rightline{$\bbox$}}

\def\ra{\rangle}

\def\sg{\sigma}

\def\sg2{\sigma^2}

\def\__{_{\infty}}

\begin{document}

\title{
A central limit theorem for weighted averages of spins
in the high temperature region of the Sherrington-Kirkpatrick model.
}

\author{Dmitry Panchenko \thanks{
Department of Mathematics,
Massachusetts Institute of Technology, Cambridge, Massachusetts.}
\thanks{
Suggested running head: Weighted averages of spins.
}
\\
77 Massachusetts Avenue, Room 2-181\\
Cambridge, MA, 02451\\
email: panchenk@math.mit.edu\\
phone: 617-253-2665\\
fax: 617-253-4358\\
}

\maketitle
\begin{abstract}
In this paper we prove that in the high temperature
region of the Sherrington-Kirkpatrick model for a typical realization
of the disorder the weighted average of spins
$\sum_{i\leq N} t_i \sigma_i$ will be approximately Gaussian
provided that $\max_{i\leq N}|t_i|/\sum_{i\leq N} t_i^2$ is small.
\end{abstract}

\vspace{0.5cm}


\section{Introduction.}

Consider a space of {\it configurations} $\Sigma_N=\{-1,+1\}^N.$ 
A configuration $\sigma\in\Sigma_N$ is a vector $(\sigma_1,\ldots,\sigma_N)$
of {\it spins} $\sigma_i$ each of which can take the values $\pm1.$
Consider an array $(g_{ij})_{i,j\leq N}$ of i.i.d. standard normal
random variables that is called the {\it disorder}.
Given parameters
$\beta>0$ and $h\geq 0,$ let us
define a Hamiltonian on $\Sigma_N$
$$
-H_{N}(\sigma)=\frac{\beta}{\sqrt{N}}
\sum_{1\leq i<j\leq N} g_{ij}\sigma_i\sigma_j
+h\sum_{i\leq N}\sigma_i,\,\,\,\,\,
\sigma=(\sigma_1,\ldots,\sigma_N)\in\Sigma_N
$$
and define a Gibbs' measure $G$ on $\Sigma_N$ by
$$
G(\{\sigma\})=\exp(-H_N(\sigma))/Z_N,\,\,\,
\mbox{ where }\,\,\, 
Z_N=\sum_{\sigma\in\Sigma_N}\exp(-H_N(\sigma)).
$$
The normalizing factor
$Z_N$ is called the {\it partition function.}
Gibbs' measure
$G$ is a random measure on $\Sigma_N$ since it depends
on the disorder $(g_{ij}).$
The parameter $\beta$ physically represents the inverse
of the temperature and
in this paper we will consider only the (very) high temperature
region of the Sherrington-Kirkpatrick model which corresponds to
\begin{equation}
\beta<\beta_0
\label{hightemp}
\end{equation}
for some small absolute constant $\beta_0>0.$
The actual value $\beta_0$ is not specified here
but, in principal, it can be determined through careful
analysis of all arguments of this paper and 
references to other papers.

For any $n\geq 1$ and
a function $f$ on the product space $(\Sigma_N^n,G^{\otimes n}),$ 
$\la f\ra$ 
will denote its expectation with respect to $G^{\otimes n}$  
$$
\la f \ra = \sum_{\Sigma_N^n} f(\sigma^1,\ldots,\sigma^n) 
G^{\otimes n}(\{(\sigma^1,\ldots,\sigma^n)\}).
$$
The Sherrington-Kirkpatrick model has been studied extensively
over the past thirty years (see, for example, [1]-[8], [11]-[18]).
In this paper we will prove the following result 
concerning the high temperature region (\ref{hightemp}).

Given a vector $(t_1,\ldots,t_N)$ such that
\begin{equation}
t_1^2+\ldots + t_N^2 =1
\label{norm}
\end{equation}
let us consider a random variable on $(\Sigma_N,G)$
defined as
\begin{equation}
X=t_1\sigma_1 +\ldots + t_N\sigma_N.
\label{X}
\end{equation}
The main goal of this paper is to show that
in the high temperature region (\ref{hightemp})
the following holds.
If $\max_{i\leq N} |t_i|$ is small then for
a typical realization of the disorder $(g_{ij})$
the random variable $X$ is approximately
Gaussian r.v. with mean $\la X \ra$ and variance
$\la X^2\ra - \la X\ra^2.$ 
By the ``typical realization'' we understand that
the statement holds on the set of measure close to $1.$

This result is the analogue of a very classical
result for independent random variables. Namely,
given a sequence of independent random variables
$\xi_1,\ldots,\xi_N$
satisfying some integrability conditions
the random variable $\xi_1+\ldots+\xi_N$
will be approximately Gaussian if 
$\max_{i\leq N} \mbox{Var}(\xi_i)/\sum_{i\leq N} \mbox{Var}(\xi_i)$
is small (see, for example, \cite{Petrov} ). 
In particular, if $\sigma_1,\ldots,\sigma_N$ in (\ref{X}) 
were i.i.d. Bernoulli random variables
then $X$ would be approximately Gaussian provided that
$\max_{i\leq N} |t_i|$ is small. 

It is important to note at this point that the main  
claim of this paper in some sense is a well expected result
since it is well known that in the high temperature region the spins
become ``decoupled'' in the limit $N\to\infty.$ 
For example, Theorem 2.4.10 in \cite{SG} states that
for a fixed $n\geq 1,$ for a typical realization of the disorder $(g_{ij})$ 
the distribution $G^{\otimes n}$ becomes a product measure
when $N\to\infty.$ 
Thus, in the very essence the claim that $X$ in (\ref{X})
is approximately Gaussian is a central limit theorem
for weakly dependent random variables.
However, the entire sequence 
$(\sigma_1,\ldots,\sigma_N)$ is a much more complicated object
than a fixed finite subset $(\sigma_1,\ldots,\sigma_n),$
and some unexpected complications arise that we will try
to describe after we state our main result - Theorem \ref{main} below.

Instead of dealing with the random variable $X$ we will look at its
symmetrized version $Y = X- X',$ where $X'$ is an independent
copy of $X.$ If we can show that $Y$ is approximately
Gaussian then, obviously, $X$ will also be approximately
Gaussian. The main reason to consider a symmetrized version
of $X$ is very simple - 
it makes it much easier to keep track of numerous indices
in all the arguments below,
even though it would be possible to carry out similar
arguments for a centered version $X -\la X \ra.$

In order to show that for a typical realization
$(g_{ij})$ and a small $\max_{i\leq N} |t_i|,$ 
$Y$ is approximately Gaussian 
with mean $0$ and variance $\la Y^2 \ra$
we will proceed by showing that its moments behave like moments of 
a Gaussian random variable, i.e. 
\begin{equation}
\la Y^l \ra \approx a(l) \la Y^2\ra^{l/2},
\label{up}
\end{equation}
where 
$a(l)=\e g^l,$ 
for a standard normal random variable $g.$
Since the moments of the standard normal random variable
are also characterized by the recursive formulas
$$
a(0)=0, a(1)=1 \mbox{ and } a(l)=(l-1)a(l-2),
$$
(\ref{up}) is equivalent to
$$
\la Y^l \ra \approx (l-1)\la Y^2\ra \la Y^{l-2}\ra.
$$

Let us define
two sequences $(\sigma^{1(l)})_{l\geq 0}$ and
$(\sigma^{2(l)})_{l\geq 0}$ of jointly independent 
random variables with Gibbs' distribution $G.$
We will assume that all indices $1(l)$ and $2(l)$
are different and one 
can think of $\sigma^{1(l)}$ and $\sigma^{2(l)}$
as different coordinates of the infinite product space
$(\Sigma_N^{\infty},G^{\otimes \infty}).$
Let us define a sequence $S_l$ by
\begin{equation}
S_l = \sum_{i=1}^N t_i \si_i^l,\,\,\,
\mbox{ where }\,\,\, \si^l=\sigma^{1(l)}-\sigma^{2(l)}.
\label{Sl}
\end{equation}
In other words, $S_l$ are independent copies of $Y.$

The following Theorem is the main result of the paper.
\begin{theorem}\label{main}
There exists $\beta_0>0$ such that 
for $\beta<\beta_0$
the following holds.
For any natural numbers $n\geq 1$ and $k_1,\ldots,k_n\geq 0$ and
$k=k_1+\ldots+ k_n,$ we have
\begin{equation}   
\bigl|
\e \la \prod_{l=1}^{n} (S_l)^{k_l} \ra - 
\prod_{l=1}^n a(k_l) 
\e
\la
S_1^2
\ra^{k/2}
\bigr|
=
\O(\max_{i\leq N} |t_i|),
\label{CLT}
\end{equation}
where $\O(\cdot)$ depends on $\beta_0, n, k$ but not
on $N.$
\end{theorem}

{\bf Remark.} Theorem \ref{main} answers the question
raised in the Research problem 2.4.11 in \cite{SG}.

Theorem \ref{main} easily implies that
\begin{equation}   
\e
\Bigl(
\la \prod_{l=1}^{n} (S_l)^{k_l} \ra - 
\prod_{l=1}^n a(k_l) 
\la
S_1^2
\ra^{k/2}
\Bigr)^2
=
\O(\max_{i\leq N} |t_i|).
\label{CLT2}
\end{equation}
Indeed,
\begin{eqnarray*}
&&
\e
\Bigl(
\la \prod_{l=1}^{n} (S_l)^{k_l} \ra - 
\prod_{l=1}^n a(k_l) 
\la
S_1^2
\ra^{k/2}
\Bigr)^2
\\
&&
=
\e
\la \prod_{l=1}^{n} (S_l)^{k_l} \ra^2 -
2\Bigl(\prod_{l=1}^n a(k_l) \Bigr)
\e
\la \prod_{l=1}^{n} (S_l)^{k_l} \ra
\la S_1^2 \ra^{k/2}
+
\Bigl(\prod_{l=1}^n a(k_l) \Bigr)^2
\e 
\la S_1^2 \ra^{k}
\end{eqnarray*}
For the first and second terms on the right hand side 
we can use independent copies to represent the powers
of $\la\cdot \ra^l$ and then apply Theorem \ref{main}.
For example,
$$
\e \la \prod_{l=1}^{n} (S_l)^{k_l} \ra^2
=
\e \la \prod_{l=1}^{n} (S_l)^{k_l} 
\prod_{l=n+1}^{2n} (S_l)^{k_{l-n}} \ra
=
\bigl(\prod_{l=1}^n a(k_l) \bigr)^2
\e 
\la S_1^2 \ra^{k}+\O(\max_{i}|t_i|).
$$
Similarly,
$$
\e
\la \prod_{l=1}^{n} (S_l)^{k_l} \ra
\la S_1^2 \ra^{k/2}
=
\bigl(\prod_{l=1}^n a(k_l) \bigr)
\e 
\la S_1^2 \ra^{k}+\O(\max_{i}|t_i|).
$$
Clearly, combining these equations proves (\ref{CLT2}).
Now one can show that 
for $N\to\infty$ and $\max_{i\leq N}|t_i|\to 0$
the characteristic 
function of $(S_1,\ldots,S_n)$ can be approximated
by the characteristic function of $n$ independent
Gaussian random variables with variance $\la S_1^2\ra,$
for $(g_{ij})$ on the set of measure converging to $1$. 
Given (\ref{CLT2}) this should be a mere exercise
and we omit the details.
This, of course, implies that $(S_1,\ldots,S_n)$
are approximately independent Gaussian random variables
with respect to the measure $G^{\otimes \infty}$ and, in particular,
$S_1=\sum_{i\leq N} t_i \si_i$ is approximately Gaussian
with respect to the measure $G^{\otimes 2}.$

Theorem \ref{main} looks very similar to
the central limit theorem for the overlap
$$
R_{1,2}=\frac{1}{N} \sum_{i=1}^{N}\sigma_i^1\sigma_i^1,
$$
where $\sigma^1,\sigma^2$ are two independent copies
of $\sigma$ (see, for example, Theorem 2.3.9 and Section 2.7
in \cite{SG}). In fact, in our proofs we
follow the main ideas and techniques of Sections 2.4 - 2.7
in \cite{SG}. However, the proof of the central limit theorem
for $X$ in (\ref{X})
turned out to be by at least an order of magnitude 
more technically involved than the proof of the central limit theorem
for the overlap $R_{1,2}$ (at least we do not know any easier proof). 
One of the main reasons why the situation here gets more complicated
is the absence of symmetry. Let us try to explain this informally.
When dealing with the overlaps $R_{i,j}$ one considers the quantity
of the following type
\begin{equation}
\e \la \prod_{i<j} R_{i,j}^{k_{i,j}} \ra
\label{overlap}
\end{equation}
and approximates it by the simpler
quantities using a kind of Taylor's expansion. 
At the second order of approximation there 
appear the terms that have ``smaller complexity''
and a term that is the factor of (\ref{overlap});
one then can solve for (\ref{overlap})
and proceed by induction on the ``complexity''.
The main reason this trick works is the symmetry.
It doesn't happen in the setting of Theorem \ref{main}
due to the lack of symmetry. Instead, we will have to consider
both terms on the left hand side of (\ref{CLT}),
\begin{equation}
\e \la \prod_{l=1}^{n} (S_l)^{k_l} \ra
\mbox{ and }
\prod_{l=1}^n a(k_l) \e\la S_1^2 \ra^{k/2},
\label{twoqua}
\end{equation}
approximate both of them by a kind of Taylor's expansion
up to the fourth order and carefully keep track of
all the terms. Surprisingly, at the forth order of approximation 
some of the terms will not be small enough to yield the claim of 
Theorem \ref{main} but the ``large '' terms corresponding
to the two quantities (\ref{twoqua}) will cancel each other.

Another difficulty that arises from the lack of symmetry
is that unlike in the case of overlaps $R_{i,j}$
we can not compute explicitly the
expectation $\la X\ra$ and variance $\la X^2 \ra - \la X\ra^2.$
Finally, we will need to develop the cavity method with two coordinates
which, loosely speaking, makes two coordinates $\sigma_i,\sigma_j$
of $\sigma$ independent of all other coordinates.
In the central limit theorem for the overlaps
$R_{i,j}$ the cavity method with one coordinate 
was sufficient.

\section{Preliminary results.}

We will first state several results from \cite{SG}
that will be constantly used throughout the paper.
Lemmas $1$ through $6$ below are either
taken directly from \cite{SG} or almost identical to some of the
results \cite{SG} and, therefore, we will state them without
the proof.

Let us consider
$$
g_t(\sigma)=\sqrt{t}\Bigl(
\sigma_N
\frac{\beta}{\sqrt{N}}\sum_{i\leq N-1}
g_{iN}\sigma_i
\Bigr) + \beta \sqrt{1-t} z\sqrt{q}\sigma_N,
$$
where $z$ is a standard normal r.v. 
independent of the disorder $(g_{ij})$
and $q$ is the unique solution of the equation
\begin{equation}
q=\e\mbox{th}^2(\beta z\sqrt{q}+h).
\label{q}
\end{equation}
For $0\leq t\leq 1$ let us consider the Hamiltonian
\begin{equation}
-H_{N,t}(\sigma)=\frac{\beta}{\sqrt{N}}
\sum_{1\leq i<j\leq N-1} g_{ij}\sigma_i\sigma_j
+ g_t(\sigma)
+h\sum_{i\leq N}\sigma_i
\label{Ht}
\end{equation}
and define Gibbs' measure $G_t$ and expectation
$\la \cdot \ra_t$ similarly to $G$ and $\la \cdot \ra$
above, only using the Hamiltonian $-H_{N,t}(\sigma).$
For any $n\geq 1$ and a function $f$ on $\Sigma_N^n$ let us define
$$
\nu_t(f)=\e \la f\ra_t .
$$
The case $t=1$ corresponds to the Hamiltonian $-H_{N}(\sigma),$
and the case $t=0$ has a very special property that
the last coordinate $\sigma_N$ is independent
of the other coordinates which is the main idea of the
{\it cavity method} (see \cite{SG}).
(Cavity method is a classical and fruitful idea
in Physics (\cite{M}), but in this paper we refer to a specific
version of the cavity method invented by Talagrand.)

Given indices $l,l',$ let us define
$$
R_{l,l'}=\frac{1}{N}\sum_{i=1}^N \sigma_i^l \sigma_i^{l'}
\,\,\,\,\mbox{ and }\,\,\,\,
R_{l,l'}^{-}=\frac{1}{N}\sum_{i=1}^{N-1} \sigma_i^l \sigma_i^{l'}.
$$
The following Lemma holds.

\begin{lemma}\label{Derivative0}
For $0\leq t < 1,$ and for all functions $f$
on $\Sigma_N^n$ we have
\begin{eqnarray}
\nu_t^{\prime}(f) 
&=& 
\beta^2\sum_{1\leq l< l' \leq n}
\nu_t(f\sigma_{N}^{l} \sigma_{N}^{l'}(R_{l,l'}^{-} - q))
-\beta^2 n \sum_{l\leq n}
\nu_t(f\sigma_{N}^{l}\sigma_{N}^{n+1}(R_{l,n+1}^{-}-q))
\nonumber
\\
&+&
\beta^2\frac{n(n+1)}{2}
\nu_t(f\sigma_{N}^{n+1}\sigma_{N}^{n+2}(R_{n+1,n+2}^{-}-q)).
\label{derivative00}
\end{eqnarray}
\end{lemma}
This is Proposition 2.4.5 in \cite{SG}.
\qed

\begin{lemma} \label{RE}
There exists $\beta_0> 0$ and $L>0$
such that for $\beta<\beta_0$ and for any 
$k\geq 1,$
\begin{equation}
\nu\Bigl((R_{1,2}-q)^{2k}\Bigr)\leq 
\Bigl(\frac{Lk}{N}\Bigr)^k\,\,\,
\mbox{ and }\,\,\, 
\nu\Bigl((R_{1,2}^{-}-q)^{2k}\Bigr)\leq 
\Bigl(\frac{Lk}{N}\Bigr)^k.
\label{Rexp}
\end{equation}
\end{lemma}
This is Theorem 2.5.1 and Lemma 2.5.2 in \cite{SG}.
\qed

Roughly speaking, this two results explain the main idea behind the
key methods of \cite{SG} - the {\it cavity method} and 
the {\it smart path method}. The Hamiltonian (\ref{Ht}) 
represents a ``smart path'' between the measures $G$ and $G_0,$ 
since along this path the derivative $\nu_{t}^{\prime}(f)$
is small, because all terms in (\ref{derivative00}) 
contain a factor $R_{l,l'} - q$ which is small due to (\ref{Rexp}).
Measure $G_0$ has a special coordinate (cavity) $\sigma_N$ that
is independent of the other coordinates, which in many cases
makes it easier to analyze $\nu_0(f).$

This two lemmas imply the following Taylor expansion
for $\nu(f).$
\begin{lemma}\label{Taylor}
For a function $f$ on $\Sigma_N^{n}$ we have
\begin{equation}
\Bigl|\nu(f)-\sum_{j=0}^{m} \frac{\nu_{0}^{(j)}(f)}{j!}\Bigr|
\leq \frac{K(n)}{N^{(m+1)/2}}\nu(f^2)^{1/2}.
\label{taylor}
\end{equation}
\end{lemma}
{\bf Proof.} Proof is almost identical to Proposition 2.5.3
in \cite{SG}.

\qed

{\bf Cavity method with two coordinates.}
In this paper we will use another case of the cavity method
with two coordinates $\sigma_N, \sigma_{N-1}$ playing the special role.  
In this new case we will consider a ``smart path'' that
makes both coordinates $\sigma_N$
and $\sigma_{N-1}$ independent of other coordinates and of each other.
This is done by slightly modifying the definition of the
Hamiltonian (\ref{Ht}).
Since it will always be clear from the context
which ``smart path'' we are using,
we will abuse the notations and use the same notations as
in the case of the Hamiltonian (\ref{Ht}).

Let us consider
\begin{eqnarray*}
g_t(\sigma)
&=&
\sqrt{t}\Bigl(
\sigma_N
\frac{\beta}{\sqrt{N}}\sum_{i\leq N-2}
g_{iN}\sigma_i
+
\sigma_{N-1}
\frac{\beta}{\sqrt{N}}\sum_{i\leq N-2}
g_{i(N-1)}\sigma_i
+\frac{\beta}{\sqrt{N}}g_{(N-1)N}\sigma_{N-1}\sigma_{N}
\Bigr) 
\\
&+& 
\beta \sqrt{1-t} 
(z_1\sqrt{q} \sigma_N + z_2 \sqrt{q}\sigma_{N-1}),
\end{eqnarray*}
where $z_1, z_2$ are standard normal r.v. 
independent of the disorder $(g_{ij}).$

For $0\leq t\leq 1$ let us now consider the Hamiltonian
\begin{equation}
-H_{N,t}(\sigma)=\frac{\beta}{\sqrt{N}}
\sum_{1\leq i<j\leq N-2} g_{ij}\sigma_i\sigma_j
+ g_t(\sigma) 
+h\sum_{i\leq N}\sigma_i
\label{Ht2}
\end{equation}
and define Gibbs' measure $G_t$ and expectation
$\la \cdot \ra_t$ similarly to $G$ and $\la \cdot \ra$
above, only using the Hamiltonian (\ref{Ht2}.)
For any $n\geq 1$ and a function $f$ on $\Sigma_N^n$ let us define
$$
\nu_t(f)=\e \la f\ra_t .
$$
We will make one distinction in the notations 
between the cases (\ref{Ht}) and (\ref{Ht2}).
Namely, for $t=0$ in the case of the Hamiltonian (\ref{Ht2}) we will
denote
\begin{equation}
\la f\ra_{00} = \la f\ra_t \Bigr|_{t=0}\,\,\,
\mbox{ and }\,\,\,\,
\nu_{00}(f) = \nu_t(f)\Bigr|_{t=0}.
\label{nu00}
\end{equation}
It is clear that with respect to the Gibbs' measure $G_0$
the last two coordinates $\sigma_N$ and $\sigma_{N-1}$
are independent of the other coordinates and of each other.

Given indices $l,l'$ let us define
$$
R_{l,l'}^{=}=\frac{1}{N}
\sum_{i\leq N-2} \sigma_{i}^l \sigma_i^{l'}.
$$
The following lemma is the analogue of Lemma \ref{Derivative0}
for the case of the Hamiltonian (\ref{Ht2}).

\begin{lemma}\label{Derivative}
Consider $\nu_t(\cdot)$ that corresponds to the Hamiltonian
(\ref{Ht2}). Then,
for $0\leq t < 1,$ and for all functions $f$
on $\Sigma_N^n$ we have
\begin{equation}
\nu_t^{\prime}(f) = \mbox{I} +\mbox{II} + \mbox{III}
\label{derivative}
\end{equation}
where
\begin{eqnarray}
\mbox{I} 
&=& 
\beta^2\sum_{1\leq l< l' \leq n}
\nu_t(f\sigma_{N}^{l} \sigma_{N}^{l'}(R_{l,l'}^{=} - q))
-\beta^2 n \sum_{l\leq n}
\nu_t(f\sigma_{N}^{l}\sigma_{N}^{n+1}(R_{l,n+1}^{=}-q))
\nonumber
\\
&+&
\beta^2\frac{n(n+1)}{2}
\nu_t(f\sigma_{N}^{n+1}\sigma_{N}^{n+2}(R_{n+1,n+2}^{=}-q)),
\label{derivative1}
\end{eqnarray}
\begin{eqnarray}
\mbox{II} 
&=& 
\beta^2\sum_{1\leq l< l' \leq n}
\nu_t(f\sigma_{N-1}^{l} \sigma_{N-1}^{l'}(R_{l,l'}^{=} - q))
-\beta^2 n \sum_{l\leq n}
\nu_t(f\sigma_{N-1}^{l}\sigma_{N-1}^{n+1}(R_{l,n+1}^{=}-q))
\nonumber
\\
&+&
\beta^2\frac{n(n+1)}{2}
\nu_t(f\sigma_{N-1}^{n+1}\sigma_{N-1}^{n+2}(R_{n+1,n+2}^{=}-q)),
\label{derivative2}
\end{eqnarray}
\begin{eqnarray}
\mbox{III} 
&=& 
\frac{1}{N}
\beta^2\sum_{1\leq l< l' \leq n}
\nu_t(f\sigma_{N}^{l} \sigma_{N}^{l'}\sigma_{N-1}^{l} \sigma_{N-1}^{l'})
-\frac{1}{N}
\beta^2 n \sum_{l\leq n}
\nu_t(f\sigma_{N}^{l}\sigma_{N}^{n+1}\sigma_{N-1}^{l}\sigma_{N-1}^{n+1})
\nonumber
\\
&+&
\frac{1}{N}
\beta^2\frac{n(n+1)}{2}
\nu_t(f\sigma_{N}^{n+1}\sigma_{N}^{n+2}\sigma_{N-1}^{n+1}\sigma_{N-1}^{n+2}).
\label{derivative3}
\end{eqnarray}
\end{lemma}
{\bf Proof.} The proof repeats the proof of Proposition 2.4.5
in \cite{SG} almost without changes.
\qed

\begin{lemma}\label{RE2}
There exists $\beta_0> 0$ and $L>0$
such that for $\beta<\beta_0$ and for any 
$k\geq 1,$
\begin{equation}
\nu_{00}\Bigl((R_{1,2}^{=}-q)^{2k}\Bigr)
\leq
L\nu\Bigl((R_{1,2}^{=}-q)^{2k}\Bigr)\leq 
\Bigl(\frac{Lk}{N}\Bigr)^k.
\label{Rexp2}
\end{equation}
\end{lemma}
The second inequality is similar to (\ref{Rexp})
and it follows easily from it since $|R_{1,2}-R_{1,2}^{=}|\leq 2/N$
(see, for example, the proof of Lemma 2.5.2 in \cite{SG}).
The first inequality follows easily from Lemma \ref{Derivative}
(see, for example, Proposition 2.4.6 in \cite{SG}).

\qed

Lemma \ref{Taylor} above also holds in the case of the
Hamiltonian (\ref{Ht2}).
\begin{lemma}\label{Taylor2}
For a function $f$ on $\Sigma_N^{n}$ we have
\begin{equation}
\Bigl|\nu(f)-\sum_{j=0}^{m} \frac{\nu_{00}^{(j)}(f)}{j!}\Bigr|
\leq \frac{K(n)}{N^{(m+1)/2}}\nu(f^2)^{1/2}.
\label{taylor2}
\end{equation}
\end{lemma}
The proof is almost identical to the proof of
Proposition 2.5.3 in \cite{SG}.

\qed

To prove Theorem \ref{main} we will need several preliminary results.

First, it will be very important to control the size of
the random variables $S_l$ and we will start 
by proving exponential integrability of $S_l.$

\begin{theorem}\label{exponent}
There exist $\beta_0>0$ and $L>0$ such that for all
$\beta\leq\beta_0,$ and for all $k\geq 1$
\begin{equation}
\nu\Bigl(\bigl(\sum_{i=1}^N t_i \si_i
\bigr)^{2k}\Bigr)\leq (Lk)^k.
\label{tight}
\end{equation}
\end{theorem}
The statement of Theorem \ref{exponent}
is, obviously, equivalent to
$$
\nu\Bigl(\exp\Bigl(L^{-1}\bigl(\sum_{i=1}^N t_i \si_i\bigr)^2
\Bigr)\Bigr)\leq L,
$$
for large enough $L.$

{\bf Proof.} 
The proof mimics the proof of Theorem 2.5.1 in \cite{SG}
(stated in Lemma \ref{RE} above).

We will prove Theorem \ref{exponent} by induction over $k.$ Our induction
assumption will be the following: there exist $\beta_0>0$
and $L>0$ such that for all $\beta\leq\beta_0,$ all $N\geq 1,$
all sequences $(t_1,\ldots,t_N)$ such that $\sum_{i=1}^N t_i^2=1$
and $0\leq l\leq k,$ we have 
\begin{equation}
\nu\Bigl(\bigl(\sum_{i=1}^N t_i \si_i
\bigr)^{2l}\Bigr)\leq (Ll)^l.
\label{ind}
\end{equation}
Let us start by proving this statement for $k=1.$
We have
$$
\nu\Bigl(\bigl(\sum_{i=1}^N t_i \si_i \bigr)^{2}\Bigr)
\leq 4\sum_{i=1}^N t_i^2 +\sum_{i\not = j}
t_i t_j \nu(\si_i\si_j)\leq
4 +(\sum_{i=1}^N t_i)^2 
\nu(\si_1\si_N)\leq
4+N\nu(\si_1\si_N).
$$
Thus we need to prove that $\nu(\si_1\si_N)\leq LN^{-1},$
for some absolute constant $L>0.$
(\ref{taylor}) implies that
$$
|\nu(\si_1\si_N)-\nu_0(\si_1\si_N)-\nu_0^{\prime}(\si_1\si_N)|
\leq
LN^{-1}.
$$
We will now show that 
$\nu_0(\si_1\si_N)=0$
and $\nu_0^{\prime}(\si_1\si_N)= \O(N^{-1}).$
The fact that $\nu_0(\si_1\si_N)=0$ is obvious since
for measure $G_0^{\otimes 2}$ the last coordinates $\sigma_N^1, \sigma_N^2$
are independent of the first $N-1$ coordinates and
$\nu_0(\si_1\si_N)=\nu_0(\si_1)\nu_0(\sigma_N^1-\sigma_N^2)=0.$
To prove that $\nu_0^{\prime}(\si_1\si_N)=\O (N^{-1})$ we use
Lemma \ref{Derivative0} which in this case implies that
\begin{eqnarray*}
\nu_0^{\prime}(\si_1\si_N)
&=&
\beta^2\nu_0(\si_1\si_N \sigma_N^1\sigma_N^2 (R_{1,2}^{-}-q))
-2\beta^2\nu_0(\si_1\si_N \sigma_N^1\sigma_N^3 (R_{1,3}^{-}-q))
\\
&-&
2\beta^2\nu_0(\si_1\si_N \sigma_N^2\sigma_N^3 (R_{2,3}^{-}-q))
+3\beta^2\nu_0(\si_1\si_N \sigma_N^3\sigma_N^4 (R_{3,4}^{-}-q))
\\
&=&
\beta^2\nu_0(\si_1(\sigma_N^2-\sigma_N^1) (R_{1,2}^{-}-q))
-2\beta^2\nu_0(\si_1(\sigma_N^3-\sigma_N^1\sigma_N^2\sigma_N^3) 
(R_{1,3}^{-}-q))
\\
&-&
2\beta^2\nu_0(\si_1(\sigma_N^1\sigma_N^2\sigma_N^3 - \sigma_N^3) 
(R_{2,3}^{-}-q))
+3\beta^2\nu_0(\si_1(\sigma_N^1 -\sigma_N^2) 
\sigma_N^3\sigma_N^4 (R_{3,4}^{-}-q))
\end{eqnarray*}
Since for a fixed disorder (r.v. $g_{ij}$ and $z$) 
the last coordinates $\sigma_N^i, i\leq 4$ are independent
of the first $N-1$ coordinates and independent of each other,
we can write
\begin{eqnarray*}
&&
\nu_0^{\prime}(\si_1\si_N)
=
\beta^2 \e \Bigl(
\la \si_1(R_{1,2}^{-} - q) \ra_0
\la \sigma_N^2 - \sigma_N^1 \ra_0
-2
\la \si_1(R_{1,3}^{-} - q) \ra_0
\la \sigma_N^3 - \sigma_N^1\sigma_N^2\sigma_N^3 \ra_0
\\
&&
-2
\la \si_1(R_{2,3}^{-} - q) \ra_0
\la \sigma_N^1\sigma_N^2\sigma_N^3 -\sigma_N^3 \ra_0
+3
\la \si_1(R_{3,4}^{-} - q) \ra_0
\la \sigma_N^1 -\sigma_N^2\ra_0
\la \sigma_N^3\sigma_N^4 \ra_0
\Bigr).
\end{eqnarray*}
First of all, the first and the last terms are equal to zero
because $\la\sigma_N^1 -\sigma_N^2 \ra_0 =0.$ Next, by symmetry
$$
\la \si_1(R_{1,3}^{-} - q) \ra_0
=
\la (\sigma_1^1 -\sigma_1^2)(R_{1,3}^{-} - q) \ra_0
=
\la (\sigma_1^2 -\sigma_1^1)(R_{2,3}^{-} - q) \ra_0
=
-\la \si_1(R_{2,3}^{-} - q) \ra_0 .
$$
Therefore, we get
\begin{eqnarray*}
\nu_0^{\prime}(\si_1\si_N)
&=&
-4\beta^2 \e \Bigl(
\la \si_1(R_{1,3}^{-} - q) \ra_0
\la \sigma_N^3 - \sigma_N^1\sigma_N^2\sigma_N^3 \ra_0
\Bigr)
\\
&=&
-4\beta^2 \nu_0 (\si_1(R_{1,3}^{-} - q) )
\nu_0(\sigma_N^3 - \sigma_N^1\sigma_N^2\sigma_N^3).
\end{eqnarray*}
It remains to show that $\nu_0 (\si_1(R_{1,3}^{-} - q) )=\O(N^{-1}).$
In order to avoid introducing new notations we notice 
that it is equivalent
to proving that $\nu (\si_1(R_{1,3} - q)) = \O(N^{-1}).$
Indeed, if we are able to prove that
\begin{equation}
\forall \beta\leq \beta_0\,\, \forall N\geq 1\,\,\,
\nu (\si_1(R_{1,3} - q)) = \O(N^{-1})
\label{nothing1}
\end{equation}
then making a change of variables $N\to N-1,$
$\beta\to \beta_{-}= \beta\sqrt{1-1/N}<\beta_0,$ and $q\to q_{-},$ where
$q_{-}$ is the solution of (\ref{q}) with $\beta$ substituted
with $\beta_{-},$ we would get
$$
\nu_0 (\si_1(R_{1,3}^{-} - q_{-}) )=\O((N-1)^{-1})=\O(N^{-1}).
$$
Lemma 2.4.15 in \cite{SG} states that for $\beta\leq\beta_0,$
$|q-q_{-}|\leq LN^{-1}$ and, therefore, the above inequality
would imply that 
$\nu_0 (\si_1(R_{1,3}^{-} - q) )=\O(N^{-1}).$
To prove (\ref{nothing1}) we notice that by symmetry
$\nu (\si_1(R_{1,3} - q))=\nu (\si_N(R_{1,3} - q)),$
and we apply (\ref{taylor})
which in this case implies that
$$
|\nu (\si_N(R_{1,3} - q)) - \nu_0 (\si_N(R_{1,3} - q))|
\leq \frac{L}{\sqrt{N}} \nu ((R_{1,3} - q)^2)^{1/2}\leq
\frac{L}{N},
$$
where in the last inequality we used (\ref{Rexp}).
Finally,
$$
\nu_0 (\si_N(R_{1,3} - q))=
-q\nu_0(\si_N)+
\frac{1}{N}\nu_0 (\si_N \sigma_N^1\sigma_N^3 )
+\nu_0 (\si_N R_{1,3}^{-})) = 
\frac{1}{N}\nu_0 (\si_N \sigma_N^1\sigma_N^3 )=
\O(N^{-1}).
$$
This finishes the proof of (\ref{ind}) for $k=1.$
It remains to prove the induction step.
One can write
\begin{equation}
\nu\Bigl(\bigl(\sum_{i=1}^N t_i \si_i
\bigr)^{2k+2}\Bigr)=
\sum_{i=1}^N t_i  
\nu\Bigl(\si_i\bigl(\sum_{j=1}^N t_j \si_j
\bigr)^{2k+1}\Bigr).
\label{sum}
\end{equation}
Let us define $\nu_i(\cdot)$ in the same way
we defined $\nu_0(\cdot)$ only now the $i$-th coordinate plays
the same role as the $N$-th coordinate played for $\nu_0.$
Using Proposition 2.4.7 in \cite{SG} we get that
for any $\tau_1,\tau_2>1$ such that $1/\tau_1 + 1/\tau_2 =1,$
\begin{eqnarray}
&&
\Bigl|\nu\Bigl(\si_i\bigl(\sum_{j=1}^N t_j \si_j
\bigr)^{2k+1}\Bigr)
-
\nu_i\Bigl(\si_i\bigl(\sum_{j=1}^N t_j \si_j
\bigr)^{2k+1}\Bigr)\Bigr|
\nonumber
\\
&&
\leq
L\beta^2
\nu\Bigl(\bigl(\sum_{i=1}^N t_i \si_i
\bigr)^{\tau_1(2k+1)}\Bigr)^{1/\tau_1}
\nu(|R_{1,2}-q|^{\tau_2})^{1/\tau_2}.
\label{est}
\end{eqnarray}
Let us take $\tau_1=(2k+2)/(2k+1)$ and $\tau_2=2k+2.$
By (\ref{Rexp}) we can estimate 
\begin{equation}
\nu\Bigl(|R_{1,2} - q|^{\tau_2}\Bigr)^{1/\tau_2}\leq
L\sqrt{\frac{\tau_2}{N}}= L\sqrt{\frac{2k+2}{N}}.
\label{exp}
\end{equation}
Next, we can write
\begin{equation}
\nu\Bigl(\bigl(\sum_{i=1}^N t_i \si_i
\bigr)^{\tau_1(2k+1)}\Bigr)^{1/\tau_1}
=
\nu\Bigl(\bigl(\sum_{i=1}^N t_i \si_i
\bigr)^{(2k+2)}\Bigr)
\nu\Bigl(\bigl(\sum_{i=1}^N t_i \si_i
\bigr)^{(2k+2)}\Bigr)^{-1/(2k+2)}.
\label{prod}
\end{equation}
If for some $\beta$ and $N$ 
$$
\nu\Bigl(\bigl(\sum_{i=1}^N t_i \si_i
\bigr)^{(2k+2)}\Bigr)\leq
(k+1)^{k+1}
$$
then for this parameters the induction step is not needed
since this inequality is precisely what we are trying to prove.
Thus, without loss of generality, we can assume that 
$\nu\Bigl(\bigl(\sum_{i=1}^N t_i \si_i
\bigr)^{(2k+2)}\Bigr)\geq
(k+1)^{k+1},$ which implies that
$$
\nu\Bigl(\bigl(\sum_{i=1}^N t_i \si_i
\bigr)^{(2k+2)}\Bigr)^{-1/(2k+2)}
\leq \frac{1}{\sqrt{k+1}}.
$$
Combining this with (\ref{est}), (\ref{exp}) and (\ref{prod})
we get
$$
\nu\Bigl(\si_i\bigl(\sum_{j=1}^N t_j \si_j
\bigr)^{2k+1}\Bigr)
\leq
\nu_i\Bigl(\si_i\bigl(\sum_{j=1}^N t_j \si_j
\bigr)^{2k+1}\Bigr)+
\frac{L\beta^2}{\sqrt{N}} 
\nu\Bigl(\bigl(\sum_{j=1}^N t_j \si_j
\bigr)^{(2k+2)}\Bigr).
$$
Plugging this estimate into (\ref{sum}) we get
\begin{eqnarray*}
\nu\Bigl(\bigl(\sum_{i=1}^N t_i \si_i
\bigr)^{2k+2}\Bigr)
&\leq&
\sum_{i=1}^N t_i  
\nu_i\Bigl(\si_i\bigl(\sum_{j=1}^N t_j \si_j
\bigr)^{2k+1}\Bigr)+
\sum_{i=1}^N t_i  
\frac{L\beta^2}{\sqrt{N}} 
\nu\Bigl(\bigl(\sum_{i=1}^N t_i \si_i
\bigr)^{(2k+2)}\Bigr)
\\
&\leq&
\sum_{i=1}^N t_i  
\nu_i\Bigl(\si_i\bigl(\sum_{j=1}^N t_j \si_j
\bigr)^{2k+1}\Bigr)+
L\beta^2 
\nu\Bigl(\bigl(\sum_{i=1}^N t_i \si_i
\bigr)^{(2k+2)}\Bigr),
\end{eqnarray*}
since (\ref{norm}) implies that $\sum t_i \leq \sqrt{N}.$
If $L\beta^2\leq 1/2,$ this implies that
\begin{equation}
\nu\Bigl(\bigl(\sum_{i=1}^N t_i \si_i
\bigr)^{2k+2}\Bigr)
\leq
2\sum_{i=1}^N t_i  
\nu_i\Bigl(\si_i\bigl(\sum_{j=1}^N t_j \si_j
\bigr)^{2k+1}\Bigr).
\label{halfind}
\end{equation}
One can write,
$$
\nu_i\Bigl(\si_i\bigl(\sum_{j=1}^N t_j \si_j
\bigr)^{2k+1}\Bigr)
=
\nu_i\Bigl(\si_i\Bigl(\bigl(\sum_{j=1}^N t_j \si_j
\bigr)^{2k+1}
-\bigl(\sum_{j\not = i} t_j \si_j \bigr)^{2k+1}
\Bigr)
\Bigr),
$$
since
$\nu_i\Bigl(\si_i
\bigl(\sum_{j\not = i} t_j \si_j \bigr)^{2k+1}
\Bigr)=0.$
Using the inequality 
$$
|x^{2k+1} - y^{2k+1}|\leq (2k+1) |x - y|(x^{2k} + y^{2k})
$$
we get
\begin{equation}
\nu_i\Bigl(\si_i\bigl(\sum_{j=1}^N t_j \si_j
\bigr)^{2k+1}\Bigr)
\leq 4(2k+1)t_i 
\Bigl[
\nu_i\Bigl(\bigl(\sum_{j\not =i} t_j \si_j
\bigr)^{2k}\Bigr)+
\nu_i\Bigl(\bigl(\sum_{j=1}^N t_j \si_j
\bigr)^{2k}\Bigr)
\Bigr].
\label{quarterind}
\end{equation}
First of all, by induction hypothesis (\ref{ind}) we have
$$
\nu_i\Bigl(\bigl(\sum_{j\not =i} t_j \si_j
\bigr)^{2k}\Bigr)\leq (Lk)^k,
$$
since this is exactly (\ref{ind}) for parameters $N-1,$
$\beta_{-}=\beta\sqrt{1-1/N},$ and since $\sum_{j\not = i}t_j^2\leq 1.$
Next, by Proposition 2.4.6 in \cite{SG} we have
$$
\nu_i\Bigl(\bigl(\sum_{j=1}^N t_j \si_j
\bigr)^{2k}\Bigr)\leq L
\nu\Bigl(\bigl(\sum_{j=1}^N t_j \si_j
\bigr)^{2k}\Bigr)\leq (Lk)^k,
$$
where in the last inequality we again used (\ref{ind}).
Thus, (\ref{halfind}) and (\ref{quarterind}) imply
$$
\nu\Bigl(\bigl(\sum_{i=1}^N t_i \si_i
\bigr)^{2k+2}\Bigr)
\leq
16(2k+1)\sum_{i=1}^N t_i^2 (Lk)^k \leq
32(k+1) (Lk)^k \leq (L(k+1))^{k+1},
$$
for $L$ large enough.
This completes the proof of the induction step and Theorem
\ref{exponent}.

\qed

{\bf Remark.}
Theorem \ref{exponent} and Lemmas \ref{RE} and \ref{RE2}
will be often used implicitly in the proof of Theorem \ref{main}
in the following way. For example, if we consider
a sequence $S_l$ defined in (\ref{Sl}) then
by H\"older's inequality (first with respect to 
$\la\cdot \ra$ and then with respect to $\e$) one can write
$$
\nu\bigl(
(R_{1,2}-q)(R_{2,3}-q) S_1^2 S_2^6
\bigr)
\leq
\nu\bigl((R_{1,2}-q)^4\bigr)^{1/4}
\nu\bigl((R_{2,3}-q)^4\bigr)^{1/4}
\nu\bigl(S_1^8\bigr)^{1/4}
\nu\bigl(S_1^{12}\bigr)^{1/4}=
\O(N^{-1}),
$$
where in the last equality we applied
Theorem \ref{exponent} and Lemma \ref{RE}.
Similarly, when we consider a function that is a product
of the factors of the type $R_{l,l'} -q$ or $S_l,$
we will simply say that each factor $R_{l,l'} - q$
contributes $\O(N^{-1/2})$ and each factor $S_l$ contributes
$\O(1).$

The following result plays the central role
in the proof of Theorem \ref{main}.
We consider a function
$$
\phi=\prod_{l=1}^n (S_l)^{q_l},
$$
where $S_l$ are defined in (\ref{Sl})
and where $q_l$ are arbitrary natural numbers,
and we consider the following quantity
$$
\nu\bigl((R_{l,l'}-q)(R_{m,m'}-q)\phi\bigr).
$$
We will show that this quantity essentially 
does not depend on the choice of pairs $(l,l')$ and $(m,m')$ or,
more accurately, it depends only on their joint configuration.
This type of quantities will appear when one considers the
second derivative of $\nu\bigl(\phi\bigr),$
after two applications of Lemma 1 or Lemma 4,
and we will be able to cancel some of these terms
up to the smaller order approximation.

\begin{lemma}\label{Lemma1}
There exists $\beta_0>0$ such that for $\beta<\beta_0$
the following holds.
Consider four pairs of indices $(l,l'), (m,m'), (p,p')$
and $(r,r')$ such that none of them is equal to 
$(1(j),2(j))$ for $j\leq n.$ Then,
if either $(l,l')\not = (m,m')$ and  $(p,p')\not = (r,r')$
or $(l,l')= (m,m')$ and $(p,p')= (r,r')$
then
\begin{equation}
\nu\bigl((R_{l,l'}-q)(R_{m,m'}-q)\phi\bigr) -
\nu\bigl((R_{p,p'}-q)(R_{r,r'}-q)\phi\bigr)
= \O(\max|t_i| N^{-1}),
\label{secondorder}
\end{equation}
where $\O(\cdot)$ depends on $n,\beta_0,\sum_{l\leq n} q_l$
but not on $N.$
\end{lemma}
{\bf Proof.}
The proof is based on the following observation.
Given $(l,l')$ consider
\begin{equation}
T_{l,l'}=N^{-1}(\sigma^l - b)\cdot (\sigma^{l'}-b),
T_l = N^{-1} (\sigma^l - b)\cdot b,
T = N^{-1} b\cdot b -q,
\label{T}
\end{equation}
where $b=\la \sigma \ra = (\la \sigma_i \ra)_{i\leq N}.$
One can express $R_{l,l'} - q$ as
\begin{equation}
R_{l,l'} - q = T_{l,l'} + T_{l} + T_{l'} + T.
\label{RasT}
\end{equation}
The joint behavior of these quantities (\ref{T}) was completely described 
in Sections 6 and 7 of \cite{SG}.
Our main observation here is that under the restrictions on indices 
made in the statement of Lemma \ref{Lemma1} the function $\phi$
will be ``almost'' independent of these quantities and all
proofs in \cite{SG} can be carried out with some minor modifications.
Let us consider the case when
$(l,l')\not = (m,m')$ and  $(p,p')\not = (r,r').$
Using (\ref{RasT}) we can write
$(R_{l,l'}-q)(R_{m,m'}-q)$ 
as the sum of terms of the following types:
$$
T_{l,l'}T_{m,m'},\, T_{l,l'}T_{m},\, T_{l,l'}T,\,
T_l T_m,\, T_l T \mbox{ and } TT.
$$ 
Similarly, we can decompose
$(R_{p,p'}-q)(R_{r,r'}-q).$
The terms on the left hand side of (\ref{secondorder})
containing a factor $TT$ will obviously cancel out.
Thus, we only need to prove that any other term multiplied by
$\phi$ will produce a quantity of order 
$\O(\max |t_i| N^{-1}).$
Let us consider, for example, the term $\nu(T_{l,l'}T_{m,m'}\phi).$
To prove that
$\nu(T_{l,l'}T_{m,m'}\phi)=\O(\max |t_i| N^{-1})$ we will
follow the proof of Proposition 2.6.5 in \cite{SG}
with some necessary adjustments.
Let us consider indices $i(1),i(2), i(3), i(4)$ that are not equal to
any of the indices that appear in $T_{l,l'},T_{m,m'}$ or $\phi.$
Then we can write,
\begin{eqnarray}
&&
\nu(T_{l,l'}T_{m,m'}\phi)=\nu\bigl(
N^{-1}(\sigma^l - \sigma^{i(1)})\cdot(\sigma^{l'} - \sigma^{i(2)})
N^{-1}(\sigma^m - \sigma^{i(3)})\cdot(\sigma^{m'} - \sigma^{i(4)})\phi\bigr)
\nonumber
\\
&&
=
N^{-1}\sum_{j=1}^N \nu\bigl(
(\sigma_j^l - \sigma_j^{i(1)})\cdot(\sigma_j^{l'} - \sigma_j^{i(2)})
(R_{m,m'} - R_{m,i(4)} - R_{m',i(3)} + R_{i(3),i(4)})\phi
\bigr).
\label{TLL}
\end{eqnarray}
Let us consider one term in this sum, for example,
\begin{equation}
\nu\bigl(
(\sigma_N^l - \sigma_N^{i(1)})\cdot(\sigma_N^{l'} - \sigma_N^{i(2)})
(R_{m,m'} - R_{m,i(4)} - R_{m',i(3)} + R_{i(3),i(4)})\phi
\bigr).
\label{oneterm}
\end{equation}
If we define
$$
R_{l,l'}^{-}=\frac{1}{N}\sum_{i=1}^{N-1} \sigma_i^l\sigma_i^{l'},\,\,\,
S_l^{-}=\sum_{i=1}^{N-1} t_i\si_i^l,
\phi^{-}=\prod_{l=1}^n (S_l^{-})^{q_l},
$$
then we can decompose (\ref{oneterm}) as
\begin{eqnarray}
&&
\nu\bigl(
(\sigma_N^l - \sigma_N^{i(1)})(\sigma_N^{l'} - \sigma_N^{i(2)})
(R_{m,m'}^{-} - R_{m,i(4)}^{-} - R_{m',i(3)}^{-} + R_{i(3),i(4)}^{-})\phi^{-}
\bigr)
\nonumber
\\
&&
+N^{-1}\nu\bigl(
(\sigma_N^l - \sigma_N^{i(1)})(\sigma_N^{l'} - \sigma_N^{i(2)})
(\sigma_N^m - \sigma_N^{i(3)})(\sigma_N^{m'} - \sigma_N^{i(4)})
\phi^{-}
\bigr)
\nonumber
\\
&&
+ t_N R_1 + t_N^2 R_2 + \O(t_N^3 + t_N N^{-1}),
\label{decompose}
\end{eqnarray}
where $R_1$ is the sum of terms of the following type
$$
R_1^j=
\nu\bigl(
(\sigma_N^l - \sigma_N^{i(1)})(\sigma_N^{l'} - \sigma_N^{i(2)})
(\sigma_N^{1(j)} - \sigma_N^{2(j)})
(R_{m,m'}^{-} - R_{m,i(4)}^{-} - R_{m',i(3)}^{-} + R_{i(3),i(4)}^{-})
\phi_{j}^{-}
\bigr),
$$
where $\phi_{j}^{-}= \prod_{l=1}^n (S_l^{-})^{q_l}/ S_{j}^{-},$
and $R_2$ is the sum of terms of the following type
\begin{eqnarray*}
&&
R_2^{j,k}=\nu\bigl(
(\sigma_N^l - \sigma_N^{i(1)})(\sigma_N^{l'} - \sigma_N^{i(2)})
(\sigma_N^{1(j)} - \sigma_N^{2(j)})
(\sigma_N^{1(k)} - \sigma_N^{2(k)})
\\
&&
\mbox{\hspace{1.5cm}}
\times
(R_{m,m'}^{-} - R_{m,i(4)}^{-} - R_{m',i(3)}^{-} + R_{i(3),i(4)}^{-})
\phi_{j,k}^{-}
\bigr),
\end{eqnarray*}
where 
$\phi_{j,k}^{-}= \prod_{l=1}^n (S_l^{-})^{q_l}/
(S_{j}^{-} S_{k}^{-}).$
First of all,
$$
|R_2^{j,k}| \leq L \nu\bigl((R_{1,2}^{-} - q)^2\bigr)^{1/2}
\nu\bigl((\phi_{j,k}^{-})^2\bigr)^{1/2}=\O(N^{-1/2}),
$$
using Theorem \ref{exponent} and Lemma \ref{RE}.
To bound $R_1^j$ we notice that $\nu_0 (R_1^j ) = 0,$
and, moreover, $\nu_0^{\prime}(R_1^j) = \O(N^{-1}),$
since by (\ref{derivative00})  
each term in the derivative will have another factor
$R_{l,l'}^{-}-q.$
Therefore, using (\ref{taylor}) we get
$$
\nu (R_1^j) = \O(N^{-1}).
$$
The second term in (\ref{decompose}) will
have order $\O(N^{-3/2})$ 
since
$$
\nu_0\bigl(
(\sigma_N^l - \sigma_N^{i(1)})(\sigma_N^{l'} - \sigma_N^{i(2)})
(\sigma_N^m - \sigma_N^{i(3)})(\sigma_N^{m'} - \sigma_N^{i(4)})
\phi^{-}
\bigr)=0
$$
and one can again apply (\ref{taylor}).
Thus the last two lines in (\ref{decompose}) will be of order
$$
\O(t_N^3 + t_N^2 N^{-1/2} + t_N N^{-1} + N^{-3/2}).
$$
To estimate the first term in (\ref{decompose})
we apply Proposition 2.6.3 in \cite{SG} which in this case implies
\begin{eqnarray*}
&&
\nu\bigl(
(\sigma_N^l - \sigma_N^{i(1)})(\sigma_N^{l'} - \sigma_N^{i(2)})
(R_{m,m'}^{-} - R_{m,i(4)}^{-} - R_{m',i(3)}^{-} + R_{i(3),i(4)}^{-})\phi^{-}
\bigr)
\\
&&
=   
L\beta^2 
\nu\bigl(
(R_{l,l'} - R_{l,i(2)} - R_{l',i(1)} + R_{i(1),i(2)})
(R_{m,m'}^{-} - R_{m,i(4)}^{-} - R_{m',i(3)}^{-} + R_{i(3),i(4)}^{-})\phi^{-}
\bigr)
\\
&&
+\O(N^{-3/2}).
\end{eqnarray*}
Now, using the similar decomposition as (\ref{oneterm}),
(\ref{decompose}) one can easily show that
\begin{eqnarray*}
&&
\nu\bigl(
(R_{l,l'} - R_{l,i(2)} - R_{l',i(1)} + R_{i(1),i(2)})
(R_{m,m'}^{-} - R_{m,i(4)}^{-} - R_{m',i(3)}^{-} + R_{i(3),i(4)}^{-})\phi^{-}
\bigr)
\\
&&
=\nu\bigl(
(R_{l,l'} - R_{l,i(2)} - R_{l',i(1)} + R_{i(1),i(2)})
(R_{m,m'} - R_{m,i(4)} - R_{m',i(3)} + R_{i(3),i(4)})\phi
\bigr)
\\
&&
+\O(N^{-3/2}+t_N N^{-1})
=\nu(T_{l,l'}T_{m,m'}\phi)+
\O(N^{-3/2}+t_N N^{-1}).
\end{eqnarray*}
Thus, combining all the estimates the term (\ref{oneterm})
becomes
\begin{eqnarray*}
&&
\nu\bigl(
(\sigma_N^l - \sigma_N^{i(1)})\cdot(\sigma_N^{l'} - \sigma_N^{i(2)})
(R_{m,m'} - R_{m,i(4)} - R_{m',i(3)} + R_{i(3),i(4)})\phi
\bigr)
\\
&&
=L\beta^2 \nu(T_{l,l'}T_{m,m'}\phi) +
\O(t_N^3 + t_N^2 N^{-1/2} + t_N N^{-1} + N^{-3/2}).
\end{eqnarray*}
All other terms on the right-hand side of (\ref{TLL})
can be written in exactly the same way, by using the cavity method in
the corresponding coordinate and, thus, (\ref{TLL}) becomes
\begin{eqnarray*}
&&
\nu(T_{l,l'}T_{m,m'}\phi)=
\sum_{j=1}^N N^{-1}
\Bigl(
L\beta^2 \nu(T_{l,l'}T_{m,m'}\phi) +
\O(t_j^3 + t_j^2 N^{-1/2} + t_j N^{-1} + N^{-3/2})
\Bigr)
\\
&&
=
L\beta^2 \nu(T_{l,l'}T_{m,m'}\phi) +
\O(\max|t_i|N^{-1}).
\end{eqnarray*}
For small enough $\beta,$ e.g. $L\beta^2\leq 1/2$
this implies that $\nu(T_{l,l'}T_{m,m'}\phi)=\O(\max |t_i|N^{-1}).$
To prove (\ref{secondorder}) in the case when
$(l,l')\not = (m,m')$ and  $(p,p')\not = (r,r'),$
it remains to estimate all other terms produces by decomposition
(\ref{RasT}) and this is done by following the proofs
of corresponding results in the Section 2.6 of \cite{SG}.

The case when $(l,l')= (m,m')$ and  $(p,p')= (r,r')$
is slightly different. The decomposition of
$(R_{l,l'}-q)^2$ using (\ref{RasT})
will produce new terms
$\nu(T_{l,l'}^2 \phi)$ and $\nu(T_l^2 \phi),$
which are not small but up to the terms of order
$\O(\max |t_i|N^{-1})$ will be equal to the corresponding terms
produces by the decomposition of $(R_{p,p'}-q)^2.$
To see this, once again, one should follow the proofs of the corresponding 
results in the Section 2.6 of \cite{SG} with minor
changes.

\qed

\section{Proof of Theorem \ref{main}}

Theorem \ref{main} is obvious if at least one $k_l$ is odd
since in this case
the left hand side of (\ref{CLT}) will be equal to $0.$
We will assume that all $k_l$ are even and, moreover,
at least one of them is greater than $2,$ say $k_1\geq 4.$
Since $a(l)=(l-1)a(l-2),$ in order to prove (\ref{CLT})
it is, obviously, enough to prove
\begin{equation}
\bigl|
\e \la \prod_{l=1}^{n} (S_l)^{k_l} \ra - 
(k_1 - 1) 
\e \la
(S_0)^2 (S_1)^{k_1 - 2} 
\prod_{l=2}^n
(S_l)^{k_l}
\ra
\bigr|
=
\O(\max_{i\leq N} |t_i|).
\label{CLTind}
\end{equation}
We will try to analyze and compare the terms on the left hand side.
Let us write
\begin{equation}
\prod_{l=1}^{n} (S_l)^{k_l}=\sum_{i=1}^{N} t_i \si_i^{1}
(S_1)^{k_1 -1}\prod_{l=2}^{n} (S_l)^{k_l}
\label{first}
\end{equation}
and 
\begin{equation}
(S_0)^2 (S_1)^{k_1 - 2} \prod_{l=2}^n (S_l)^{k_l}=
\sum_{i=1}^{N} t_i \si_i^{0}
(S_0) (S_1)^{k_1 - 2} \prod_{l=2}^n (S_l)^{k_l}.
\label{second}
\end{equation}
From now on we will carefully analyze terms in (\ref{first})
in several steps and at each step we will notice that one of two
things happens: 

(a) The term produced at the same step of our analysis
carried out for (\ref{second})
is exactly the same up to a constant $k_1 - 1$;

(b) The term is ``small'' meaning that after combining all the steps
one would get something of order $\O(\max |t_i|).$

Obviously these observations will imply (\ref{CLTind}).

Let us look at one term in (\ref{first}) and (\ref{second}),
for example,
\begin{equation}
\si_N^{1}(S_1)^{k_1 -1}\prod_{l=2}^{n} (S_l)^{k_l}
\,\,\,
\mbox{ and }\,\,\,
\si_N^{0}
(S_0) (S_1)^{k_1 - 2} \prod_{l=2}^n (S_l)^{k_l}.
\label{lastterms}
\end{equation}
If we define $S_l^{-}$ by the equation
$$
S_l = S_l^{-} + t_N \si_N^l,
$$
then,
\begin{eqnarray*}
&&
\si_N^{1}(S_1)^{k_1 -1}\prod_{l=2}^{n} (S_l)^{k_l}
=\si_N^{1}(S_1^{-}+t_N \si_N^1)^{k_1 -1}\prod_{l=2}^{n} 
(S_l^{-} + t_N \si_N^l)^{k_l}
\\
&&
=
\si_N^{1}(S_1^{-})^{k_1 -1}\prod_{l=2}^{n} (S_l^{-})^{k_l}
+(k_1 -1) t_N (\si_N^1)^2 (S_1^{-})^{k_1 -2}\prod_{l=2}^{n} (S_l^{-})^{k_l}
\\
&&
+
t_N\sum_{l=2}^n k_l \si_N^1\si_N^l (S_1^{-})^{k_1 -1}
(S_l^{-})^{k_l -1}\prod_{j\not = 1,l} (S_j^{-})^{k_j} 
+\O(t_N^2)
\\
&&
= \mbox{I} +t_N \mbox{II} +t_N \mbox{III} +\O(t_N^2).
\end{eqnarray*}
and
\begin{eqnarray*}
&&
\si_N^{0} (S_0)(S_1)^{k_1 -2}\prod_{l=2}^{n} (S_l)^{k_l}
=\si_N^{0}(S_0^{-} + t_N \si_N^{0})
(S_1^{-}+t_N \si_N^1)^{k_1 -2}\prod_{l=2}^{n} 
(S_l^{-} + t_N \si_N^l)^{k_l}
\\
&&
=
\si_N^{0}(S_0^{-})(S_1^{-})^{k_1 -2}\prod_{l=2}^{n} (S_l^{-})^{k_l}
+ t_N (\si_N^0)^2 (S_1^{-})^{k_1 -2}\prod_{l=2}^{n} (S_l^{-})^{k_l}
\\
&&
+
t_N 
\Bigl(
\sum_{l=2}^n k_l \si_N^0 \si_N^l (S_1^{-})^{k_1 -2}
(S_l^{-})^{k_l -1}\prod_{j\not = 1,l} (S_j^{-})^{k_j} 
\\
&&
+
(k_1 - 2) \si_N^0 \si_N^1 (S_1^{-})^{k_1 -3}
\prod_{j\not = 2}^n (S_j^{-})^{k_j} 
\Bigr)
+\O(t_N^2)
\\
&&
= \mbox{IV} +t_N \mbox{V} +t_N \mbox{VI} +\O(t_N^2).
\end{eqnarray*}
First of all, $\nu_0 (\mbox{III}) = \nu_0(\mbox{VI}) = 0$  and, therefore,
applying (\ref{taylor})
$$
t_N \nu(\mbox{III}) = \O(t_N N^{-1/2})\,\,\, \mbox{ and }\,\,\,
\nu(\mbox{VI}) = \O(t_N N^{-1/2}).
$$
Next, again using (\ref{taylor})
\begin{eqnarray*}
t_N \nu(\mbox{II}) 
&=& 
t_N \nu_0(\mbox{II}) + t_N\O(N^{-1/2})
\\
&=&
t_N (k_1 -1) \nu_0((\si_N^1)^2)\nu_0
((S_1^{-})^{k_1 -2}\prod_{l=2}^{n} (S_l^{-})^{k_l})
+\O(t_N N^{-1/2})
\end{eqnarray*}
and
\begin{eqnarray*}
t_N \nu(\mbox{V}) 
&=& 
t_N \nu_0(\mbox{V}) + t_N\O(N^{-1/2})
\\
&=&
t_N \nu_0((\si_N^1)^2)\nu_0
((S_1^{-})^{k_1 -2}\prod_{l=2}^{n} (S_l^{-})^{k_l})
+\O(t_N N^{-1/2}).
\end{eqnarray*}
Thus the contribution of the terms II and V in (\ref{CLTind}) will cancel out
- the first appearance of case (a) mentioned above.

The terms of order $\O(t_N^2 + t_N N^{-1/2})$ 
when plugged back into (\ref{first}) and (\ref{second}) 
will produce 
\begin{equation}
\sum_{i=1}^N t_i \O(t_i^2 + t_i N^{-1/2}) =
\O(\max_{i\leq N} |t_i| + N^{-1/2})
= \O(\max_{i\leq N} |t_i|).
\label{effect1}
\end{equation}
Here we, of course, assume that similar analysis is carried out
for the $i$-th term in (\ref{first}) and (\ref{second})
with the only difference that the $i$th coordinate plays the special role 
in the definition of $\nu_0.$

We now proceed to analyze the terms I and IV.
If we define $S_l^{=}$ by the equation
$$
S_l^{-} = S_l^{=} + t_{N-1} \si_{N-1}^l,
$$
then,
\begin{eqnarray*}
&&
\mbox{I}=
\si_N^{1}(S_1^{-})^{k_1 -1}\prod_{l=2}^{n} (S_l^{-})^{k_l}
=
\si_N^{1}(S_1^{=})^{k_1 -1}\prod_{l=2}^{n} (S_l^{=})^{k_l}
\\
&&
+
t_{N-1}(k_1 -1) R_0
+
t_{N-1} R_1
+t_{N-1}^2 (R_{21} + R_{22} + R_{23})
+t_{N-1}^3 R_3 + \O(t_{N-1}^4). 
\end{eqnarray*}
where
$$
R_0 = \si_N^{1}\si_{N-1}^1
(S_1^{=})^{k_1 -2}\prod_{l=2}^{n} (S_l^{=})^{k_l}
$$
$$
R_1 = \sum_{l=2}^n k_l
\si_N^1 \si_{N-1}^l (S_1^{=})^{k_1 -1}
(S_l^{=})^{k_l - 1} \prod_{j\not = 1,l} (S_l^{=})^{k_l}, 
$$
$$
R_{21} = {k_1 -1 \choose 2} 
\si_N^1 (\si_{N-1}^1)^2 (S_1^{=})^{k_1-3}
\prod_{l=1}^n (S_l^{=})^{k_l},
$$
$$
R_{22} = \sum_{l=2}^n {k_l \choose 2}
\si_N^1 (\si_{N-1}^l)^2 (S_1^{=})^{k_1-1}
(S_l^{=})^{k_l -2}
\prod_{j\not =1,l}^n (S_j^{=})^{k_l},
$$
and where $R_{23}$ is the sum of terms of the following type
$$
\si_N^1 \si_{N-1}^l \si_{N-1}^{l'}
\prod_{j=1}^n (S_j^{=})^{q_j},\,\,\,
1\leq l\not = l' \leq n
$$
for some (not important here) powers $q_l,$ and where
$R_3$ is the sum of terms of the following type
$$
\si_N^1 \si_{N-1}^l \si_{N-1}^{l'} \si_{N-1}^{l''}
\prod_{j=1}^n (S_j^{=})^{q_j},\,\,\,
1\leq l,l',l'' \leq n.
$$
Similarly,
\begin{eqnarray*}
&&
\mbox{IV}=
\si_N^{0}(S_0^{-})
(S_1^{-})^{k_1 -2}\prod_{l=2}^{n} (S_l^{-})^{k_l}
=
\si_N^{0}(S_0^{=})
(S_1^{=})^{k_1 -2}\prod_{l=2}^{n} (S_l^{=})^{k_l}
\\
&&
+
t_{N-1} \bar{R}_0
+
t_{N-1} \bar{R}_1
+t_{N-1}^2 (\bar{R}_{21} + \bar{R}_{22} + \bar{R}_{23})
+t_{N-1}^3 \bar{R}_3 + \O(t_{N-1}^4). 
\end{eqnarray*}
where
$$
\bar{R}_0 = \si_N^{0}\si_{N-1}^0
(S_1^{=})^{k_1 -2}\prod_{l=2}^{n} (S_l^{=})^{k_l}
$$
$$
\bar{R}_1 = 
\si_N^0 \si_{N-1}^1 S_0^{=}
(S_1^{=})^{k_1 -3}
\prod_{j\not = 0,1} (S_l^{=})^{k_l}
+
\sum_{l=2}^n k_l
\si_N^0 \si_{N-1}^l S_0^{=}
(S_1^{=})^{k_1 -2}
(S_l^{=})^{k_l - 1} \prod_{j\not = 0,1,l} (S_l^{=})^{k_l}, 
$$
$$
\bar{R}_{21} = 
{k_1 - 2 \choose 2}
\si_N^0 (\si_{N-1}^1)^2 (S_0^{=})
(S_1^{=})^{k_1-4}
\prod_{l=2}^n (S_l^{=})^{k_l},
$$
$$
\bar{R}_{22} = \sum_{l=2}^n {k_l \choose 2}
\si_N^0 (\si_{N-1}^l)^2 (S_0^{=})
(S_1^{=})^{k_1-2}
(S_l^{=})^{k_l -2}
\prod_{j\not =0,1,l}^n (S_j^{=})^{k_l},
$$
and where $\bar{R}_{23}$ is the sum of terms of the following type
$$
\si_N^0 \si_{N-1}^l \si_{N-1}^{l'}
\prod_{j=0}^n (S_j^{=})^{q_j},\,\,\,
0\leq l\not = l' \leq n
$$
for some (not important here) powers $q_l,$ and where
$\bar{R}_3$ is the sum of terms of the following type
$$
\si_N^0 \si_{N-1}^l \si_{N-1}^{l'} \si_{N-1}^{l''}
\prod_{j=0}^n (S_j^{=})^{q_j},\,\,\,
0\leq l,l',l'' \leq n.
$$

(Step 1). First of all since $\nu_0(R_3)=0,$ we have 
$\nu(R_3)=\O(N^{-1/2})$ and $t_{N-1}^3\nu(R_3)=\O(t_{N-1}^3 N^{-1/2}).$
Next let us show that $\nu(R_{23})=\O(N^{-1}).$
Indeed, one need to note that $\nu_{00}(R_{23})=0,$
and using Lemma \ref{Derivative}, 
$\nu_{00}^{\prime}(R_{23})=\O(N^{-1})$
since each term produced by (\ref{derivative1}) will have
a factor $\la \si_{N-1}^l\ra_{00} = 0,$ each term produced by
(\ref{derivative2}) will have a factor $\la \si_N^1 \ra_{00}=0,$
and each term produced by (\ref{derivative3}) has factor $N^{-1}.$
Thus it remains to use (\ref{taylor2})
to show that $\nu(R_{23})=\O(N^{-1}).$

Similarly, one can show that
$\nu(\bar{R}_3)=\O(N^{-1/2})$ and $\nu(\bar{R}_{23})=\O(N^{-1}).$

(Step 2).
Let us show now that $\nu(R_{1})=\O(N^{-3/2}).$
Let us consider one individual term
$$
R_{1l}=
\si_N^1 \si_{N-1}^l (T_1^{=})^{k_1 -1}
(T_l^{=})^{k_l - 1} \prod_{j\not = 1,l} (T_l^{=})^{k_l}. 
$$
Obviously, $\nu_{00}(R_{1l})=0.$ To show that 
$\nu_{00}^{\prime}(R_{1l})=0,$ let us first note that
the terms produced by (\ref{derivative1}) 
will contain a factor $\la \si_{N-1}^l\ra_{00}=0,$
the terms produced by (\ref{derivative2}) 
will contain a factor $\la \si_{N}^1\ra_{00}=0,$
and the terms produced by (\ref{derivative3}) 
will contain a factor $\la (S_1^{=})^{k_1-1}\ra_{00}=0,$
since $k_1 -1$ is odd and $S_1^{=}$ is symmetric.
For the second derivative we will have different types of terms
produced by a combination of (\ref{derivative1}),
(\ref{derivative2}) and (\ref{derivative3}).
The terms produced by using (\ref{derivative3}) twice
will have order $\O(N^{-2});$ the terms produced by
using (\ref{derivative3}) and either (\ref{derivative2})
or (\ref{derivative1}) will have order $\O(N^{-3/2}),$
since the factor $R_{l,l'}^{=}-q$ will produce $N^{-1/2};$
the terms produced by (\ref{derivative1}) and (\ref{derivative1}),
or by (\ref{derivative2}) and (\ref{derivative2}) will be equal to
$0$ since they will contain factors $\la \si_{N-1}^l \ra_{00}=0$
and $\la \si_{N}^1 \ra_{00}=0$ correspondingly.
Finally, let us consider the terms produced by 
(\ref{derivative1}) and (\ref{derivative2}), e.g.
$$
\nu_{00}\bigl(
R_{1l}
\sigma_{N}^{m}\sigma_{N}^{m'}(R_{m,m'}^{=}-q)
\sigma_{N-1}^{p}\sigma_{N-1}^{p'}(R_{p,p'}^{=}-q)
\bigr).
$$
It will obviously be equal to $0$ unless $m,p\in \{1(1),2(1)\}$
and $m',p'\in\{1(l),2(l)\}$ since, otherwise,
there will be a factor $\la \si_N^1\ra_{00}=0$ 
or $\la \si_N^l\ra_{00}=0.$
All non zero terms will cancel due to the
following observation. Consider, for example, the term
\begin{eqnarray*}
&&
\nu_{00}\bigl(
R_{1l}
\sigma_{N}^{1(1)}\sigma_{N}^{1(l)}(R_{1(1),1(l)}^{=}-q)
\sigma_{N-1}^{2(1)}\sigma_{N-1}^{2(l)}(R_{2(1),2(l)}^{=}-q)
\bigr)
=
\nu_{00}\bigl(
\sigma_{N}^{1(l)} - \sigma_{N}^{1(1)}
\sigma_{N}^{2(1)}\sigma_{N}^{1(l)}
\Bigl)\times
\\
&&
\nu_{00}\bigl(
\sigma_{N-1}^{1(1)}\sigma_{N-1}^{1(2)}
\sigma_{N-1}^{2(l)}-\sigma_{N}^{2(l)}
\Bigl)
\nu_{00}\bigr(
(R_{1(1),1(l)}^{=}-q) (R_{2(1),2(l)}^{=}-q)
(T_1^{=})^{k_1 -1}
(T_l^{=})^{k_l - 1} \prod_{j\not = 1,l} (T_l^{=})^{k_l}
\bigr),
\end{eqnarray*}
which corresponds to $m=1(1), m'=1(l),p=2(1)$ and $p'=2(l).$
There will also be a similar term that corresponds to
$m=2(1), m'=1(l),p=1(1)$ and $p'=2(l)$ (indices
$m$ and $p$ are changed)
\begin{eqnarray*}
&&
\nu_{00}\bigl(
R_{1l}
\sigma_{N}^{2(1)}\sigma_{N}^{1(l)}(R_{2(1),1(l)}^{=}-q)
\sigma_{N-1}^{1(1)}\sigma_{N-1}^{2(l)}(R_{1(1),2(l)}^{=}-q)
\bigr)
=
\nu_{00}\bigl( 
\sigma_{N}^{1(1)}\sigma_{N}^{2(1)}\sigma_{N}^{1(l)}
-\sigma_{N}^{1(l)}
\Bigl)\times
\\
&&
\nu_{00}\bigl(
\sigma_{N}^{2(l)}-
\sigma_{N-1}^{1(1)}\sigma_{N-1}^{1(2)}\sigma_{N-1}^{2(l)}
\Bigl)
\nu_{00}\bigr(
(R_{2(1),1(l)}^{=}-q) (R_{1(1),2(l)}^{=}-q)
(T_1^{=})^{k_1 -1}
(T_l^{=})^{k_l - 1} \prod_{j\not = 1,l} (T_l^{=})^{k_l}
\bigr).
\end{eqnarray*}
These two terms will cancel since the product of the first two factors
is unchanged and, making the change of variables
$1(1)\to 2(1),$ $2(1)\to 1(1)$ in the last factor we get
(note that $T_1^{=}\to - T_{1}^{=}$)
\begin{eqnarray*}
&&
\nu_{00}\bigr(
(R_{2(1),1(l)}^{=}-q) (R_{1(1),2(l)}^{=}-q)
(T_1^{=})^{k_1 -1}
(T_l^{=})^{k_l - 1} \prod_{j\not = 1,l} (T_l^{=})^{k_l}
\bigr)
\\
&&
=
\nu_{00}\bigr(
(R_{1(1),1(l)}^{=}-q) (R_{2(1),2(l)}^{=}-q)
(-T_1^{=})^{k_1 -1}
(T_l^{=})^{k_l - 1} \prod_{j\not = 1,l} (T_l^{=})^{k_l}
\bigr)
\\
&&
=
-\nu_{00}\bigr(
(R_{1(1),1(l)}^{=}-q) (R_{2(1),2(l)}^{=}-q)
(T_1^{=})^{k_1 -1}
(T_l^{=})^{k_l - 1} \prod_{j\not = 1,l} (T_l^{=})^{k_l}
\bigr).
\end{eqnarray*}
Using (\ref{taylor2}) we finally get that $\nu(R_1)=\O(N^{-3/2}).$

Similarly, one can show that $\nu(\bar{R}_1)=\O(N^{-3/2}).$

(Step 3).
Next, we will show that
\begin{equation}
\nu(R_{21}) - (k_1 -1)\nu(\bar{R}_{21})=\O(N^{-1})
\label{R21}
\end{equation}
and
\begin{equation}
\nu(R_{22}) - (k_1 -1)\nu(\bar{R}_{22})=\O(N^{-1}).
\label{R22}
\end{equation}
We will prove only (\ref{R21}) since (\ref{R22}) is proved
similarly.
Since $\nu_{00}(R_{21})=\nu_{00}(\bar{R}_{21})=0$ it is enough to
prove that
\begin{eqnarray}
&&
{k_1 -1 \choose 2} 
\nu_{00}^{\prime}\bigl(
\si_N^1 (\si_{N-1}^1)^2 (S_1^{=})^{k_1-3}
\prod_{l=2}^n (S_l^{=})^{k_l}\bigr)
\nonumber
\\
&&
=
(k_1 -1)
{k_1 - 2 \choose 2}
\nu_{00}^{\prime}\bigl(
\si_N^0 (\si_{N-1}^1)^2 (S_0^{=})
(S_1^{=})^{k_1-4}
\prod_{l=2}^n (S_l^{=})^{k_l}\bigr)
+
\O(N^{-1}).
\label{compare21}
\end{eqnarray}
On both sides the terms produced by (\ref{derivative2})
will be equal to $0,$ the terms produced by (\ref{derivative3})
will be of order $\O(N^{-1}),$ thus, it suffices to compare the terms
produced by (\ref{derivative1}).
For the left hand side the terms produced by (\ref{derivative1})
will be of the type
$$
\nu_{00}\bigl(
\si_N^1 \sigma_{N}^m \sigma_{N}^{m'}
(\si_{N-1}^1)^2 
(R_{m,m'}^{=}-q)
(S_1^{=})^{k_1-3}
\prod_{l=2}^n (S_l^{=})^{k_l}\bigr)
$$
and will be equal to $0$ unless $m\in\{1(1),2(1)\}$
and $m'\not\in \{1(1),2(1)\}.$
For a fixed $m'$ consider the sum of two terms that correspond to
$m=1(1)$ and $m=2(1),$ i.e.
\begin{eqnarray*}
&&
\nu_{00}\bigl(
(\sigma_{N}^{m'} - \sigma_{N}^{1(1)}\sigma_N^{2(1)}\sigma_N^{m'})
(\si_{N-1}^1)^2 
(R_{1(1),m'}^{=}-q)
(S_1^{=})^{k_1-3}
\prod_{l=2}^n (S_l^{=})^{k_l}\bigr)
\\
&&
+
\nu_{00}\bigl(
(\sigma_{N}^{1(1)}\sigma_N^{2(1)}\sigma_N^{m'} - \sigma_N^{m'})
(\si_{N-1}^1)^2 
(R_{2(1),m'}^{=}-q)
(S_1^{=})^{k_1-3}
\prod_{l=2}^n (S_l^{=})^{k_l}\bigr)
\\
&&
=
\nu_{00}\bigl(
(\sigma_{N}^{m'} - \sigma_{N}^{1(1)}\sigma_N^{2(1)}\sigma_N^{m'})
\bigr)
\nu_{00}\bigl(
(\si_{N-1}^1)^2 
\bigr)
\nu_{00}\bigl(
(R_{1(1),m'}^{=}- R_{2(1),m'}^{=})
(S_1^{=})^{k_1-3}
\prod_{l=2}^n (S_l^{=})^{k_l}\bigr)
\\
&&
=c
\nu_{00}\bigl(
(R_{1(1),m'}^{=}- R_{2(1),m'}^{=})
(S_1^{=})^{k_1-3}
\prod_{l=2}^n (S_l^{=})^{k_l}\bigr)
.
\end{eqnarray*}
For $m'\in\{1(2),2(2),\ldots,1(n),2(n)\}$
this term will have a factor $\beta^2,$
and for $m'=2n+1$ it will have a factor 
$-\beta^2(2n).$
Similarly, the derivative on the right hand side of
(\ref{compare21}) will consist of the terms of type
$$
c
\nu_{00}\bigl(
(R_{1(0),m'}^{=}- R_{2(0),m'}^{=})
(S_0^{=})
(S_1^{=})^{k_1-4}
\prod_{l=2}^n (S_l^{=})^{k_l}\bigr).
$$
For $m'\in\{1(1),2(1),\ldots,1(n),2(n)\}$
this term will have a factor $\beta^2,$
and for $m'=2n+3$ it will have a factor 
$-\beta^2(2n+2).$
We will show next that for any
$m'$ and $m'',$
\begin{eqnarray}
&&
\nu_{00}\bigl(
(R_{1(1),m'}^{=}- R_{2(1),m'}^{=})
(S_1^{=})^{k_1-3}
\prod_{l=2}^n (S_l^{=})^{k_l}\bigr)
\nonumber
\\
&&
=
(k_1 -3)
\nu_{00}\bigl(
(R_{1(0),m''}^{=}- R_{2(0),m''}^{=})
(S_0^{=})
(S_1^{=})^{k_1-4}
\prod_{l=2}^n (S_l^{=})^{k_l}\bigr)
+\O(N^{-1}).
\label{compare211}
\end{eqnarray}
This implies, for example, that all terms in the derivatives
are "almost" independent of the index $m'.$
This will also imply (\ref{compare21}) since,
given arbitrary fixed $m',$
the left hand side of (\ref{compare21})
will be equal to
$$
(k_1-3){k_1 -1 \choose 2}
c\beta^2 
\bigl((2n-2)-(2n)\bigr)
\nu_{00}\bigl(
(R_{1(0),m'}^{=}- R_{2(0),m'}^{=})
(S_0^{=})
(S_1^{=})^{k_1-4}
\prod_{l=2}^n (S_l^{=})^{k_l}\bigr)
+\O(N^{-1})
$$
and the right hand side of (\ref{compare21})
will be equal to
$$
(k_1-1){k_1 -2 \choose 2}
c\beta^2 
\bigl((2n)-(2n+2)\bigr)
\nu_{00}\bigl(
(R_{1(0),m'}^{=}- R_{2(0),m'}^{=})
(S_0^{=})
(S_1^{=})^{k_1-4}
\prod_{l=2}^n (S_l^{=})^{k_l}\bigr)
+\O(N^{-1}),
$$
which is the same up to the terms of order $\O(N^{-1}).$
For simplicity of notations, instead of proving (\ref{compare211})
we will prove
\begin{eqnarray}
&&
\nu\bigl(
(R_{1(1),m'}- R_{2(1),m'})
(S_1)^{k_1-3}
\prod_{l=2}^n (S_l)^{k_l}\bigr)
\nonumber
\\
&&
=
(k_1 -3)
\nu\bigl(
(R_{1(0),m''}- R_{2(0),m''})
(S_0)
(S_1)^{k_1-4}
\prod_{l=2}^n (S_l)^{k_l}\bigr)
+\O(N^{-1}).
\label{compare212}
\end{eqnarray}
Let us write the left hand side as
$$
\nu\bigl(
(R_{1(1),m'}- R_{2(1),m'})
(S_1)^{k_1-3}
\prod_{l=2}^n (S_l)^{k_l}\bigr)
=N^{-1}\sum_{i=1}^N \nu(U_i),
$$
where
$$
U_i = 
(\sigma_i^{1(1)}-\sigma_i^{2(1)})\sigma_i^{m'}
(S_1)^{k_1-3}
\prod_{l=2}^n (S_l)^{k_l}
=
\si_i^{1}\sigma_i^{m'}
(S_1)^{k_1-3}
\prod_{l=2}^n (S_l)^{k_l}
.
$$
and consider one term in this sum, for example,
$\nu(U_N).$
Using (\ref{taylor}), one can write
$$
\nu(U_N)=\nu_0(U_N)+\nu_0^{\prime}(U_N)+\O(N^{-1})
$$
and 
$$
\nu^{\prime}(U_N)=\nu_0^{\prime}(U_N)+\O(N^{-1}),
$$
since each term in the derivative already contains a factor
$R_{l,l'}^{-}-q.$
Thus,
$$
\nu(U_N)=\nu_0(U_N)+\nu^{\prime}(U_N)+\O(N^{-1}).
$$
Similarly,
$$
\nu(U_i)=\nu_i(U_i)+\nu^{\prime}(U_i)+\O(N^{-1}),
$$
where $\nu_i$ is defined the same way as $\nu_0$
only now $i$th coordinated plays the same role as 
$N$th coordinate plays for $\nu_0 (=\nu_N).$
Therefore,
\begin{eqnarray*}
&&
\nu\bigl(
(R_{1(1),m'}- R_{2(1),m'})
(S_1)^{k_1-3}
\prod_{l=2}^n (S_l)^{k_l}\bigr)
=
N^{-1}\sum_{i=1}^N \nu_i(U_i)
\\
&&
+
\nu^{\prime}\bigl(
(R_{1(1),m'}- R_{2(1),m'})
(S_1)^{k_1-3}
\prod_{l=2}^n (S_l)^{k_l}\bigr)
+\O(N^{-1})
=
N^{-1}\sum_{i=1}^N \nu_i(U_i)
+\O(N^{-1}),
\end{eqnarray*}
again using (\ref{taylor2}) and (\ref{Rexp2})
and writing
$R_{1(1),m'}- R_{2(1),m'}=
(R_{1(1),m'}-q) - (R_{2(1),m'} - q).$
Similarly one can write,
$$
\nu\bigl(
(R_{1(0),m''}- R_{2(0),m''})
(S_0)
(S_1)^{k_1-4}
\prod_{l=2}^n (S_l)^{k_l}\bigr)
=
N^{-1}\sum_{i=1}^N \nu_i(V_i) + \O(N^{-1}),
$$
where
$$
V_i=
(\sigma_i^{1(0)}-\sigma_i^{2(0)})\sigma_i^{m''}
(S_0)
(S_1)^{k_1-4}
\prod_{l=2}^n (S_l)^{k_l}.
$$
If we can finally show that
$$
\nu_i(U_i)=(k_1 -3)\nu_i(V_i)+ \O(t_i^2), 
$$
this will prove (\ref{compare212}) and (\ref{compare21}).
For example, if we consider $\nu_0(U_N),$ 
\begin{eqnarray*}
&&
\nu_0(U_N)=
\nu_0\bigl(
\si_N^{1}\sigma_N^{m'}
(S_1)^{k_1-3}\prod_{l=2}^n (S_l)^{k_l}
\Bigr)
=
\nu_0\bigl(
\si_N^{1}\sigma_N^{m'}
(S_1^{-}+t_N\si_N^1)^{k_1-3}\prod_{l=2}^n (S_l^{-}+t_N\si_N^l)^{k_l}
\Bigr)
\\
&&
=
\nu_0(\si_N^{1})
\nu_0(\sigma_N^{m'})
\nu_0\bigl(
(S_1^{-})^{k_1-3}\prod_{l=2}^n (S_l^{-})^{k_l}
\Bigr)
\\
&&
+(k_1- 3)t_N
\nu_0(
(\si_N^{1})^2)
\nu_0(\sigma_N^{m'})
\nu_0\bigl(
(S_1^{-})^{k_1-4}\prod_{l=2}^n (S_l^{-})^{k_l}
\Bigr)
\\
&&
+
t_N\sum_{l=2}^n
\nu_0(\si_N^{1}\si_N^l)
\nu_0(\sigma_N^{m'})
\nu_0\bigl(
(S_1^{-})^{k_1-3}(S_l^{-})^{k_l-1}\prod_{j\not =1,l}^n (S_j^{-})^{k_l}
\Bigr)
+\O(t_N^2)
\\
&&
=
(k_1- 3)t_N
\nu_0(
(\si_N^{1})^2)
\nu_0(\sigma_N^{m'})
\nu_0\bigl(
(S_1^{-})^{k_1-4}\prod_{l=2}^n (S_l^{-})^{k_l}
\Bigr)
+\O(t_N^2),
\end{eqnarray*}
since all other terms are equal to $0.$
Similarly, one can easily see that
$$
\nu_0(V_N)=
t_N\nu_0(
(\si_N^{1})^2)
\nu_0(\sigma_N^{m'})
\nu_0\bigl(
(S_1^{-})^{k_1-4}\prod_{l=2}^n (S_l^{-})^{k_l}
\Bigr)
+\O(t_N^2).
$$
This finishes the proof of (\ref{compare21}).

The comparison of $R_{22}$ and $\bar{R}_{22}$ can be carried out
exactly the same way.

(Step 4).
The last thing we need to prove is that
\begin{equation}
\nu(R_0)-\nu(\bar{R}_0)=\O(\max_i |t_i| N^{-1})
\label{R0dif}
\end{equation}
or, in other words,
$$
\nu\bigl(
\si_N^{1}\si_{N-1}^1
(S_1^{=})^{k_1 -2}\prod_{l=2}^{n} (S_l^{=})^{k_l}
\bigr)-
\nu\bigl(
\si_N^{0}\si_{N-1}^0
(S_1^{=})^{k_1 -2}\prod_{l=2}^{n} (S_l^{=})^{k_l}
\bigr)
=\O(\max_{i}|t_i| N^{-1}).
$$
First of all, clearly, $\nu_{00}(R_0)=\nu_{00}(\bar{R}_0)=0.$
Next we will show that 
\begin{equation}
\nu_{00}^{\prime}(R_0)-\nu_{00}^{\prime}(\bar{R}_0)=0.
\label{R1dif}
\end{equation}
The terms produced by (\ref{derivative1}) and (\ref{derivative2})
will be equal to $0,$ because they will contain either the factor
$\la \si_N^1 \ra_{00}=0$ ($\la \si_N^0 \ra_{00}=0$ ) 
or the factor $\la \si_{N-1}^1 \ra_{00}=0$ 
($\la \si_{N-1}^0 \ra_{00}=0$).
The terms of $\nu_{00}^{\prime}(R_0)$ produced
by (\ref{derivative3}) will be of the type
$$
N^{-1}\nu_{00}(
\si_N^{1} \sigma_N^{m}\sigma_N^{m'} )
\nu_{00}(\si_{N-1}^1 \sigma_{N-1}^m\sigma_{N-1}^{m'})
\nu_{00}\bigl(
(S_1^{=})^{k_1 -2}\prod_{l=2}^{n} (S_l^{=})^{k_l}
\bigr)
$$
and they will be different from $0$ only if
$m\in\{1(1),2(1)\}$ and $m'\not\in\{1(1),2(2)\}.$
For $m\in\{1(1),2(1)\}$ and $m'\in\{1(2),2(2),\ldots,1(n),2(n)\}$
these terms will have a factor $\beta^2,$
and for $m\in\{1(1),2(1)\}$ and $m'=2n+1$
these terms will have a factor $-(2n)\beta^2.$
Similarly,
the terms of $\nu_{00}^{\prime}(\bar{R}_0)$ produced
by (\ref{derivative3}) will be of the type
$$
N^{-1}\nu_{00}(
\si_N^{0} \sigma_N^{p}\sigma_N^{p'} )
\nu_{00}(\si_{N-1}^0 \sigma_{N-1}^p\sigma_{N-1}^{p'})
\nu_{00}\bigl(
(S_1^{=})^{k_1 -2}\prod_{l=2}^{n} (S_l^{=})^{k_l}
\bigr)
$$
and they will be different from $0$ only if
$p\in\{1(0),2(0)\}$ and $p'\not\in\{1(0),2(0)\}.$
For $p\in\{1(0),2(0)\}$ and $p'\in\{1(1),2(1),\ldots,1(n),2(n)\}$
these terms will have a factor $\beta^2,$
and for $p\in\{1(0),2(0)\}$ and $p'=2n+3$
these terms will have a factor $-(2n+2)\beta^2.$
For $m=1(1)$ (or $m=2(1)$) and a corresponding
$p=1(0)$ (or $p=2(0)$) the non zero terms 
above will be equal, so when
we add up the factors over $m'$ and $p'$ we get
$$
\beta^2((2n-2)-2n - (2n)+(2n+2))=0.
$$
This shows that
$\nu_{00}^{\prime}(R_0)-\nu_{00}^{\prime}(\bar{R}_0)=0.$

Next we will show that
\begin{equation}
\nu_{00}^{\prime\prime}(R_0)-\nu_{00}^{\prime\prime}(\bar{R}_0)=
\O(\max_{i} |t_i| N^{-1}).
\label{Rdif}
\end{equation}
The second derivative will have different types of terms
produced by an iterated application 
of (\ref{derivative1}),
(\ref{derivative2}) and (\ref{derivative3}).
The terms produced by using (\ref{derivative3}) twice
will have order $\O(N^{-2});$ the terms produced by
using (\ref{derivative3}) and either (\ref{derivative2})
or (\ref{derivative1}) will have order $\O(N^{-3/2}),$
since the factor $R_{l,l'}^{=}-q$ will contribute $N^{-1/2}$

via the application of (\ref{Rexp2});
the terms produced by (\ref{derivative1}) and (\ref{derivative1}),
or by (\ref{derivative2}) and (\ref{derivative2}) will be equal to
$0$ since they will contain a factor $\la \si_{N-1}^1 \ra_{00}=0$
or $\la \si_{N}^1 \ra_{00}=0$ correspondingly.
Finally, let us consider the terms produced by 
(\ref{derivative1}) and (\ref{derivative2}). For 
$\nu_{00}^{\prime\prime}(R_0)$ they will be of the type
\begin{eqnarray*}
&&
\nu_{00}\bigl(
R_{0}
\sigma_{N}^{m}\sigma_{N}^{m'}(R_{m,m'}^{=}-q)
\sigma_{N-1}^{p}\sigma_{N-1}^{p'}(R_{p,p'}^{=}-q)
\bigr)
\\
&&
=
\nu_{00}(
\si_N^{1} \sigma_N^{m}\sigma_N^{m'} )
\nu_{00}(\si_{N-1}^1 \sigma_{N-1}^p\sigma_{N-1}^{p'})
\nu_{00}\bigl(
(R_{m,m'}^{=}-q)(R_{p,p'}^{=}-q)
(S_1^{=})^{k_1 -2}\prod_{l=2}^{n} (S_l^{=})^{k_l}
\bigr)
\end{eqnarray*}
and will be equal to $0$ unless $m,p\in\{1(1),2(1)\}$
and $m',p'\not\in\{1(1),2(1)\}.$
For $\nu_{00}^{\prime\prime}(\bar{R}_0)$ the terms will be of the type
\begin{eqnarray*}
&&
\nu_{00}\bigl(
\bar{R}_{0}
\sigma_{N}^{m}\sigma_{N}^{m'}(R_{m,m'}^{=}-q)
\sigma_{N-1}^{p}\sigma_{N-1}^{p'}(R_{p,p'}^{=}-q)
\bigr)
\\
&&
=
\nu_{00}(
\si_N^{0} \sigma_N^{m}\sigma_N^{m'} )
\nu_{00}(\si_{N-1}^0 \sigma_{N-1}^p\sigma_{N-1}^{p'})
\nu_{00}\bigl(
(R_{m,m'}^{=}-q)(R_{p,p'}^{=}-q)
(S_1^{=})^{k_1 -2}\prod_{l=2}^{n} (S_l^{=})^{k_l}
\bigr)
\end{eqnarray*}
and will be equal to $0$ unless $m,p\in\{1(0),2(0)\}$
and $m',p'\not\in\{1(0),2(0)\}.$
Now, to show (\ref{Rdif})
one only needs to apply (\ref{secondorder})
and notice that for each case in Lemma \ref{Lemma1}
(i.e. for $(m,m')=(p,p')$ or $(m,m')\not =(p,p')$)
there will be equal number of positive and negative terms
that will cancel each other out up to the terms of
order $\O(\max_{i}|t_i|N^{-1}).$
The count of this terms is done similarly to what
we did in the proof of (\ref{R1dif}) and we omit it.
Finally, (\ref{R1dif}) and (\ref{Rdif}) imply (\ref{R0dif})
via the application of (\ref{taylor2}).

Now we can combine Steps $1$ through $4$ to get that
\begin{eqnarray*}
&&
\nu(\mbox{I})-(k_1 -1)\nu(\mbox{IV})=
\nu\bigl(
\si_N^{1}
(S_1^{=})^{k_1 -1}\prod_{l=2}^{n} (S_l^{=})^{k_l}
\bigr)
-(k_1 -1)
\nu\bigl(
\si_N^{0}(S_0^{=})
(S_1^{=})^{k_1 -2}\prod_{l=2}^{n} (S_l^{=})^{k_l}
\bigr)
\\
&&
+
\O(t_{N-1}^4 + t_{N-1}^3 N^{-1/2} +t_{N-1}^2 N^{-1} +
t_{N-1} \max_{i}|t_i| N^{-1}).
\end{eqnarray*}
We notice that the first two terms on the right hand side
$$
\nu\bigl(
\si_N^{1}
(S_1^{=})^{k_1 -1}\prod_{l=2}^{n} (S_l^{=})^{k_l}
\bigr)
-(k_1 -1)
\nu\bigl(
\si_N^{0}(S_0^{=})
(S_1^{=})^{k_1 -2}\prod_{l=2}^{n} (S_l^{=})^{k_l}
\bigr)
$$
are absolutely similar to 
$\nu(\mbox{I})-(k_1 -1)\nu(\mbox{IV}),$
with the only difference that $S_l^{-}$ is
now $S_l^{=}.$ Thus we can proceed by induction
to show that
\begin{eqnarray*}
&&
\nu(\mbox{I})-(k_1 -1)\nu(\mbox{IV})=
\sum_{j=1}^{N-1}
\O(t_{j}^4 + t_{j}^3 N^{-1/2} +t_{j}^2 N^{-1} +
t_{j} \max_{i}|t_i| N^{-1}).
\end{eqnarray*}
We can now add up the contributions of the terms 
$\mbox{I}$ and $\mbox{IV}$
(and terms similar to (\ref{lastterms}) 
arising from (\ref{first}) and (\ref{second}))
in the left hand side of
(\ref{CLTind}) to get
$$
\sum_{i\leq N} \sum_{j\not = i} t_i
\O(t_{j}^4 + t_{j}^3 N^{-1/2} +t_{j}^2 N^{-1} +
t_{j} \max_{l}|t_l| N^{-1})
=
\O(\max_{l}|t_l|),
$$
which is a simple calculus exercise, provided that 
$\sum_{i\leq N} t_i^2 =1.$
This, together with (\ref{effect1}),
completes the proof of (\ref{CLTind}) 
and the proof of Theorem \ref{main}.

\qed

\end{document}